# RELATIONS BETWEEN THE STRONG GLOBAL DIMENSION, COMPLEXES OF FIXED SIZE AND DERIVED CATEGORY

YOHNY CALDERÓN-HENAO, FELIPE GALLEGO-OLAYA, AND HERNÁN GIRALDO

ABSTRACT. Let $\mathbb{Z}$ be the integer numbers, $\Bbbk$ an algebraically closed field, $\Lambda$ a finite dimensional $\Bbbk$-algebra, mod$\Lambda$ the category of finitely generated right modules, proj$\Lambda$ the full subcategory of mod$\Lambda$ consisting of all projective $\Lambda$-modules, and $C_n(\text{proj}\Lambda)$ the bounded complexes of projective $\Lambda$-modules of fixed size for any integer $n \geq 2$. We found an algorithm to calculate the strong global dimension of $\Lambda$, when $\Lambda$ has finite strong global dimension and is derived-discrete, using the Auslander-Reiten quivers of the categories $C_n(\text{proj}\Lambda)$. Also, we show the relation between the Auslander-Reiten quiver of the bounded derived category $D^b(\Lambda)$ and the Auslander-Reiten quiver of $C_{\eta+1}(\text{proj}\Lambda)$, where $\eta = \text{s.gl.dim}(\Lambda)$ (strong global dimension of $\Lambda$).

## 1. INTRODUCTION

The theory of almost split sequences (or Auslander-Reiten sequences) was extended with success to the study of the bounded derived category of mod$\Lambda$, denoted by $D^b(\Lambda)$, where $\Lambda$ is a finite dimensional $\Bbbk$-algebra and $\Bbbk$ is an algebraically closed field. The notion of almost split sequences gave origin to the notion of almost split triangles (or Auslander-Reiten triangles). In his fundamental book [18], D. Happel proved that the derived category of a finite global dimension algebra has almost split triangles; he was very successful in describing the Auslander-Reiten quiver in the case of a hereditary algebra.

Let mod$\Lambda$ be the category of finitely generated right $\Lambda$-modules and let proj$\Lambda$ be the full subcategory of mod$\Lambda$ consisting of all projective $\Lambda$-modules. The category of all complexes of finitely generated projective right $\Lambda$-modules of fixed size, denoted by $C_n(\text{proj}\Lambda)$ for any integer $n \geq 2$ was introduced and studied by R. Bautista, M.J. Souto Salorio, and R. Zuazua (see [6]). In [12], C. Chaio, I. Pratti, and M.J. Souto Salorio continued with the study of this category and showed that the Auslander-Reiten quiver of $C_n(\text{proj}\Lambda)$ can be constructed by the well-known knitting algorithm used to build the Auslander-Reiten quiver mod$\Lambda$. Furthermore, in [3], R. Bautista studied the representation type of $C_n(\text{proj}\Lambda)$ for all positive integer $n \geq 2$. Recently in [11], C. Chaio, A.G. Chaio, and I. Pratti proved that if the strong global dimension of $\Lambda$ is finite ($\eta = \text{s.gl.dim}(\Lambda)$ with $\eta \in \mathbb{N}$), then, $C_{\eta+1}(\text{proj}\Lambda)$ is of finite representation type if and only if $C_n(\text{proj}\Lambda)$ is of finite representation type for each $n \geq 2$.

At the moment, there is not a general method to calculate the strong global dimension of $\Lambda$, but some works have been carried out in this direction, see [2] and [11]. In the first one, E.R. Alvares, P. Le Meur, and E.N. Marcos used the structure of the derived category to characterize the strong global dimension of $\Lambda$, when $\Lambda$ is piecewise hereditary algebra. In the second one, C. Chaio, A.G. Chaio, and I. Pratti determined the strong global dimension of some piecewise hereditary finite dimensional algebras considering their ordinary quivers with relations.

In this paper, we developed an algorithm to calculate the strong global dimension of $\Lambda$ when $\Lambda$ has finite strong global dimension and is derived-discrete. For this, we used the techniques given in [12], by C. Chaio, I. Pratti, and M.J. Souto Salorio to construct the Auslander-Reiten quiver of the category $C_n(\text{proj}\Lambda)$ for $n \geq 2$ (see Theorem 3.22, and Algorithm 3.24). Also, we proved that if $\Lambda$ is an algebra with finite strong global dimension, then, $C_{\eta+1}(\text{proj}\Lambda)$ is of tame representation type if and only if $C_n(\text{proj}\Lambda)$ is of tame representation type for all $n \geq 1$ (see Theorem 4.4); furthermore, we showed that, if $\Lambda$ has finite strong global dimension, then, up to translation, all irreducible morphisms and all almost split $\mathcal{E}_{\eta+1}$-sequences in $C_{\eta+1}(\text{proj}\Lambda)$ are irreducible morphisms and almost split $\mathcal{E}_{\eta+i}$-sequences in $C_{\eta+i}(\text{proj}\Lambda)$ for $i \geq 2$ (see Theorem 3.21). And we also proved that, up to isomorphism, some irreducible morphisms and some almost split $\mathcal{E}_{\eta+1}$-sequences in $C_{\eta+1}(\text{proj}\Lambda)$ are irreducible morphisms and Auslander-Reiten triangles in $D^b(\Lambda)$ (see Theorem 4.1). Finally, we obtained the Auslander-Reiten quiver of $D^b(\Lambda)$ as $\mathbb{Z}$ copies of a special subquiver (see Definition 4.2) of the Auslander-Reiten of $C_{\eta+1}(\text{proj}\Lambda)$ (see Theorem 4.3).

This paper is organized as follows: In Section 2, we recall some notations and preliminary results. In Section 3, we show how the indecomposable complexes shape, the irreducible morphisms, and the almost split $\mathcal{E}_{\eta+i}$-sequences in $C_{\eta+i}(\text{proj}\Lambda)$ for $i \geq 1$ are, if $\Lambda$ has finite strong global dimension (see [15] for more details). We provide an algorithm to calculate the strong global dimension of a finite dimension $\Bbbk$-algebra. In Section 4, we show the relation in terms of irreducible morphisms and almost split $\mathcal{E}_{\eta+1}$-sequences in $C_{\eta+1}(\text{proj}\Lambda)$ and Auslander-Reiten triangles $D^b(\Lambda)$ and two examples for applying the algorithm and calculating the strong global dimension and the Auslander-Reiten quiver of the bounded derived category $D^b(\Lambda)$ from the Auslander-Reiten quiver of $C_{\eta+1}(\text{proj}\Lambda)$.

## 2. Preliminary results

We started giving notations, definitions, and basic facts, which we will use in the subsequent sections. $\mathcal{A}$ will always denote an additive category, which is Krull-Schmidt and $\text{Hom}_\mathcal{A}(X, Y)$ the set of the morphism from $X$ to $Y$, where $X$ and $Y$ are objects in $\mathcal{A}$, when no confusion arises, we will use the notation $\text{Hom}(X, Y)$. If $X = Y$, $\text{Hom}(X, X) := \text{End}(X)$ and $\text{Aut}(X)$ are elements in $\text{End}(X)$, which are isomorphisms. A morphism $f\colon X \to Y$ in $\text{Hom}_\mathcal{A}(X, Y)$ is said to be a **section** (resp. **retraction**) if there is a morphism $h\colon Y \to X$ such that $hf = 1_X$ (resp. $fh = 1_Y$). When one of these conditions holds, $f$ is said to be a split morphism. A morphism $f$ is called **irreducible** when it is not split; but, for any factorization $f = hg$, $h$ is a retraction or $g$ is a section. A morphism $f\colon X \to Y$ in $\text{Hom}(X, Y)$ is called **right almost split** if $f$ is not a retraction; and if $g\colon Z \to Y$ is not a retraction, there is a $s\colon Z \to X$ with $fs = g$. A morphism $f\colon X \to Y$ in $\text{Hom}(X, Y)$ is called **minimal right** if $fu = f$ with $u \in \text{End}(X)$ implies $u \in \text{Aut}(X)$. A morphism $f\colon X \to Y$ in $\mathcal{A}$ is called **minimal right almost split** if $f$ is right almost split and $f$ is right minimal. Dually, a **minimal left almost split morphism** is defined. Be aware that $f$ is an irreducible morphism in $\mathcal{A}$, provided $f$ is minimal right (or left) almost split morphism in $\mathcal{A}$.

A pair $(i, d)$ of composable morphism $X \xrightarrow{i} Y \xrightarrow{d} Z$ in $\mathcal{A}$ is said to be **exact** if $(X, i)$ is a kernel of $d$ and $(Z, d)$ is a cokernel of $i$. Let $\mathcal{E}$ be a class of exact pairs closed under isomorphisms. If $(i, d) \in \mathcal{E}$, morphisms $i$ and $d$ are called **inflation** and **deflation** respectively, and the exact pairs in $\mathcal{E}$ are called **conflations**. The sequence $X \xrightarrow{i} Y \xrightarrow{d} Z$ in $\mathcal{A}$ is called **$\mathcal{E}$-sequence** if $(i,d) \in \mathcal{E}$. An $\mathcal{E}$-sequence $X \xrightarrow{i} Y \xrightarrow{d} Z$ is said to be **almost split** if $i$ is a minimal left almost split morphism and $d$ is a minimal right almost split morphism; in [[13], 2.2], P. Dräxler, I. Reiten, S.O. Smalø, and Ø. Solberg, show that $X$ and $Z$ are indecomposable objects in $\mathcal{A}$.

Let's recall that a complex $X = (X^i, d_X^i)_{i \in \mathbb{Z}}$ is a family of objects $X^i$ in $\mathcal{A}$ and morphisms $d_X^i \colon X^i \to X^{i+1}$ in $\mathcal{A}$, $i \in \mathbb{Z}$, such that $d_X^{i+1} d_X^i = 0$ for all $i \in \mathbb{Z}$. Morphisms $d_X^i \colon X^i \to X^{i+1}$ are called differential maps of the complex $X$. If $X$ and $Y$ are complexes, a morphism of complexes $f = (f^i)_{i \in \mathbb{Z}} \colon X \to Y$ is given by a family of morphisms $f^i \colon X^i \to Y^i$, $i \in \mathbb{Z}$, such that $d_Y^i f^i = f^{i+1} d_X^i$ for all $i \in \mathbb{Z}$. C($\mathcal{A}$) is the category whose objects are the complexes over $\mathcal{A}$, morphisms are the morphisms between two complexes, each $X^i$ is called a **cell** and each $f^i$ is called **component** of $f$ for each $i \in \mathbb{Z}$. Let C$^{-,b}(\mathcal{A})$ be the full subcategory of C($\mathcal{A}$) where the complexes are bounded above and bounded cohomology. By C$^b(\mathcal{A})$ the full subcategory of C($\mathcal{A}$) given by the bounded complexes, K($\mathcal{A}$), K$^{-,b}(\mathcal{A})$ and K$^b(\mathcal{A})$ are called homotopy categories of C($\mathcal{A}$), C$^{-,b}(\mathcal{A})$ and C$^b(\mathcal{A})$ respectively, and D$^b(\mathcal{A})$ is the bounded derived category of $\mathcal{A}$. If $X$ is a complex and $p \in \mathbb{Z}$, we have a new $X[p]$ complex given by the shifting or translation functor. For a morphism $f \colon X \to Y$ in C($\mathcal{A}$), the notation C$_f$ for the cone of the morphism $f$ is used.

In [6], R. Bautista, M.J. Souto Salorio, and R. Zuazua considered C($\mathcal{A}$) and $\mathcal{E}$ the class of exact pairs closed under isomorphisms of composable morphisms $X \xrightarrow{f} Y \xrightarrow{g} Z$ such that for all $n \in \mathbb{Z}$, the sequence $0 \to X^n \xrightarrow{f^n} Y^n \xrightarrow{g^n} Z^n \to 0$ is split exact. It is known that (C($\mathcal{A}$), $\mathcal{E}$) is an exact category in the sense of D. Quillen (see [20]), or equivalent in the sense of P. Gabriel and A.V. Roiter (see [14]), (see for instance Example 3.2 in [18]). In [[6], Notation 3.1], for $m, n \in \mathbb{Z}$ with $m < n$, C$_{[m,n]}(\mathcal{A})$ is the full subcategory of C($\mathcal{A}$) whose objects are those $X \in $ C($\mathcal{A}$) such that $X^i = 0$ if $i \notin \{m, m+1, \ldots, n\}$. This category C$_{[m,n]}(\mathcal{A})$ is a subcategory of C($\mathcal{A}$) closed under extensions and (C$_{[m,n]}(\mathcal{A})$, $\mathcal{E}_{[m,n]}$) is an exact category with $\mathcal{E}_{[m,n]}$ the class of composable morphisms in C$_{[m,n]}(\mathcal{A})$, which are in $\mathcal{E}$.

In particular, we use the notation C$_n(\mathcal{A}) := $ C$_{[1,n]}(\mathcal{A})$ and $\mathcal{E}_n = \mathcal{E}_{[1,n]}$. Also in [[6], Theorem 4.5], R. Bautista, M.J. Souto Salorio, and R. Zuazua proved the existence of the almost split $\mathcal{E}_n$-sequences in C$_n(\mathcal{A})$ or $\mathcal{E}_n$-almost split sequences. They consider $P \in \mathcal{A}$ indecomposable and defined the following three complexes in C$_n(\mathcal{A})$.

J$_k(P) = (J^i, d^i)_{i \in \mathbb{Z}}$ for $k \in \{1, \ldots, n-1\}$, with $J^i = 0$ if $i \neq k$ and $i \neq k+1$, $J^k = J^{k+1} = P$ and $d^k = 1_P$ and $d^i = 0$ in another case.

S$(P) = (X^i, d^i)_{i \in \mathbb{Z}}$ with $X^i = 0$ for $i \neq 1$, $X^1 = P$, $d^i = 0$ for all $i \in \mathbb{Z}$.

T$(P) = (X^i, d^i)_{i \in \mathbb{Z}}$ with $X^i = 0$ for $i \neq n$, $X^n = P$, $d^i = 0$ for all $i \in \mathbb{Z}$.

Following R. Bautista in [4], a complex $P$ in C$_n(\mathcal{A})$ is called $\mathcal{E}_n$-**projective** if for any deflation $f \colon Y \to Z$ and any morphism $u \colon P \to Z$ in C$_n(\mathcal{A})$ there exists $v \colon P \to Y$ such that $fv = u$. Dually, $\mathcal{E}_n$-**injective** is defined. In [[6], Corollary 3.8 and Corollary 3.9], R. Bautista, M.J. Souto Salorio, and R. Zuazua proved the following proposition.

**Proposition 2.1.** *The indecomposable complexes $\mathcal{E}_n$-projective (resp $\mathcal{E}_n$-injective) in C$_n(\mathcal{A})$ are the complexex J$_k(P)$ and T$(P)$ (resp J$_k(P)$ and S$(P)$) for $k \in \{1, \ldots, n-1\}$ with $P \in \mathcal{A}$ indecomposable.*

Let $\Bbbk$ be an algebraically closed field and $\Lambda$ a finite dimensional $\Bbbk$-algebra. By mod$\Lambda$ and proj$\Lambda$, we denote the category of finitely generated right $\Lambda$-modules and the full subcategory of finitely generated projective right $\Lambda$-modules, respectively.

We recall from D. Happel and D. Zacharia in [19], the following: If $X \in \mathrm{K}^b(\mathrm{proj}\Lambda)$ is a complex, we may consider a preimage $\bar{X}$ of $X$ in $\mathrm{C}^b(\mathrm{proj}\Lambda)$ without indecomposable projective direct summands. $\bar{X}$ is uniquely determined by $X$ up to isomorphism of bounded complexes in $\mathrm{C}^b(\mathrm{proj}\Lambda)$. Thus the following is well-defined. If $0 \neq X \in \mathrm{K}^b(\mathrm{proj}\Lambda)$ there exists $r \leq s$ such that $\bar{X}^r \neq 0 \neq \bar{X}^s$ and $\bar{X}^i = 0$ for $i < r$ and $i > s$. Then by definition the length of $X$ is defined as $\ell(X) = s - r$. We defined the **strong global dimension** of $\Lambda$ by

$$\mathrm{s.gl.dim}(\Lambda) := \sup\{\ell(X) | X \in \mathrm{K}^b(\mathrm{proj}\Lambda) \text{ indecomposable}\}.$$

For here on, $\eta$ will denote the finite strong global dimension of $\Lambda$, that is, $\eta = \mathrm{s.gl.dim}(\Lambda)$ with $\eta \geq 1$.

We know from D. Happel in [18], that $\Lambda$ is **piecewise hereditary**, if there exists a hereditary abelian category $\mathcal{H}$ such that, $\mathrm{D}^b(\Lambda)$ and $\mathrm{D}^b(\mathcal{H})$ are equivalent as triangulated categories.

The following theorem was proved in [[19], Theorem 3.2], by D. Happel and D. Zacharia.

**Theorem 2.2.** *Let $\Lambda$ be a finite dimensional $\Bbbk$-algebra. Then $\Lambda$ is a piecewise hereditary if and only if $\mathrm{s.gl.dim}(\Lambda) < \infty$.*

In this work, we focused our attention on $\mathrm{C}_n(\mathrm{proj}\Lambda)$ category, whose objects are the complexes $X$ with $X^i = 0$ if $i \notin \{1, \ldots, n\}$, i.e.,

$$X: \cdots \to 0 \to X^1 \to X^2 \to \cdots \to X^{n-1} \to X^n \to 0 \to \cdots.$$

When no confusion arises, this complex will be written as follows: $X: (X^1 \to X^2 \to \cdots \to X^{n-1} \to X^n)$, and $X^1$ is the first cell, $X^2$ is the second cell, etc and $X$ will be determined by these cells. And if $X, Y \in \mathrm{C}_n(\mathrm{proj}\Lambda)$ and $f: X \to Y$ is a morphism in $\mathrm{C}_n(\mathrm{proj}\Lambda)$, then $f$ is denoted by $f: (f^1, \ldots, f^n)$ and $f$ will be determined by these components.

We consider the following functors for $n \geq 1$.

i. $i_{-1}^n: \mathrm{C}_n(\mathrm{proj}\Lambda) \to \mathrm{C}_{n+1}(\mathrm{proj}\Lambda)$ given by $i_{-1}^n(X)^1 = 0$, $i_{-1}^n(X)^j = X^{j-1}$ for all $j \in \{2, \ldots, n+1\}$, and $d^1_{i_{-1}^n(X)} = 0$, $d^j_{i_{-1}^n(X)} = d_X^{j-1}$ for all $j \in \{2, \ldots, n\}$. If $X, Y$, and $f: X \to Y$ are in $\mathrm{C}_n(\mathrm{proj}\Lambda)$, then the morphism $i_{-1}^n(f): i_{-1}^n(X) \to i_{-1}^n(Y)$ in $\mathrm{C}_{n+1}(\mathrm{proj}\Lambda)$ is given by $i_{-1}^n(f)^1 = 0$ and $i_{-1}^n(f)^j = f^{j-1}$ for all $j \in \{2, \ldots, n+1\}$.

ii. $i_0^n: \mathrm{C}_n(\mathrm{proj}\Lambda) \to \mathrm{C}_{n+1}(\mathrm{proj}\Lambda)$ given by $i_0^n(X)^j = X^j$ for all $j \in \{1, \ldots, n\}$, $i_0^n(X)^{n+1} = 0$, and $d^j_{i_0^n(X)} = d_X^j$ for all $j \in \{1, \ldots, n-1\}$. If $X, Y$, and $f: X \to Y$ are in $\mathrm{C}_n(\mathrm{proj}\Lambda)$, then the morphism $i_0^n(f): i_0^n(X) \to i_0^n(Y)$ in $\mathrm{C}_{n+1}(\mathrm{proj}\Lambda)$ is given by $i_0^n(f)^j = f^j$ for all $j \in \{1, \ldots, n\}$.

iii. $\pi_1^n: \mathrm{C}_{n+1}(\mathrm{proj}\Lambda) \to \mathrm{C}_n(\mathrm{proj}\Lambda)$ given by $\pi_1^n(X)^j = X^{j+1}$ for all $j \in \{1, \ldots, n\}$, $d^j_{\pi_1^n(X)} = d_X^{j+1}$ for all $j \in \{1, \ldots, n-1\}$. If $X, Y$, and $f: X \to Y$ are in $\mathrm{C}_{n+1}(\mathrm{proj}\Lambda)$, then the morphism $\pi_1^n(f): \pi_1^n(X) \to \pi_1^n(Y)$ is given by $\pi_1^n(f)^j = f^{j+1}$ for all $j \in \{1, \ldots, n\}$.

iv. $\pi_0^n: \mathrm{C}_{n+1}(\mathrm{proj}\Lambda) \to \mathrm{C}_n(\mathrm{proj}\Lambda)$ given by $\pi_0^n(X)^j = X^j$ for all $j \in \{1, \ldots, n\}$, $d^j_{\pi_0^n(X)} = d_X^j$ for all $j \in \{1, \ldots, n-1\}$. If $X, Y$, and $f: X \to Y$ are in $\mathrm{C}_{n+1}(\mathrm{proj}\Lambda)$, then the morphism $\pi_0^n(f): \pi_0^n(X) \to \pi_0^n(Y)$ is given by $\pi_0^n(f)^j = f^j$ for all $j \in \{1, \ldots, n\}$.

**Remark 2.3.** We note the following facts for $n \geq 1$.

  i. $i_{-1}^n$ is translation functor $[-1]$, up to the sign of the differential maps; $i_0^n$ is inclusion functor, that is, $i_0^n = [0]$; $\pi_1^n$ is translation functor $[1]$, up to the sign of differential maps; and $\pi_0^n$ is the functor that vanish $X^{n+1}$.

  ii. $\pi_1^n$ is a left inverse functor of $i_{-1}^n$ and $\pi_0^n$ is a left inverse functor of $i_0^n$.

  iii. Let $C_{n+1}^{(1)}(\mathrm{proj}\Lambda)$ be the full subcategory of $C_{n+1}(\mathrm{proj}\Lambda)$ such that $X^1 = 0$ for all $X \in C_{n+1}(\mathrm{proj}\Lambda)$, so $\pi_1^n$ restrict to $C_{n+1}^{(1)}(\mathrm{proj}\Lambda)$ is an isomorphism with inverse isomorphism $i_{-1}^n$.

  iv. Let $C_{n+1}^{(n+1)}(\mathrm{proj}\Lambda)$ be the full subcategory of $C_{n+1}(\mathrm{proj}\Lambda)$ such that $X^{n+1} = 0$ for all $X \in C_{n+1}(\mathrm{proj}\Lambda)$, so $\pi_0^n$ restrict to $C_{n+1}^{(n+1)}(\mathrm{proj}\Lambda)$ is an isomorphism with inverse isomorphism $i_0^n$; moreover, we note that in this case $\pi_0^n = [0]$.

**Lemma 2.4.** *Let $\Lambda$ be a finite dimensional $\Bbbk$-algebra and $X$ and $Y$ complexes in $C_{n+i}(\mathrm{proj}\Lambda)$ with $i \geq 1$. If $f \colon X \to Y$ is an irreducible morphism in $C_{n+i}(\mathrm{proj}\Lambda)$, then we have the following:*

  i. *If $X^1 = 0 = Y^1$, then $\pi_1^{n+i}(f) \colon \pi_1^{n+i}(X) \to \pi_1^{n+i}(Y)$ is an irreducible morphism in $C_{n+i-1}(\mathrm{proj}\Lambda)$.*

  ii. *If $X^{n+i} = 0 = Y^{n+i}$, then $\pi_0^{n+i}(f) \colon \pi_0^{n+i}(X) \to \pi_0^{n+i}(Y)$ is an irreducible morphism in $C_{n+i-1}(\mathrm{proj}\Lambda)$.*

**Proof.** It is an easy verification. $\square$

C. Chaio, A.G. Chaio, and I. Pratti give the following definition that can be found in [[10], Definition 4.1], more exactly:

**Definition 2.5.** *Let $X = (X^i, d_X^i)_{i \in \mathbb{Z}}$ be an indecomposable complex in $C_n(\mathrm{proj}\Lambda)$.*

  i. *We say that $X$ can be **extended to the left** in $C_{n+1}(\mathrm{proj}\Lambda)$ if there is a projective $\Lambda$-module $X'^0$ and a non-zero morphism $d'^0 \colon X'^0 \to X^1$ such that $d_X^1 d'^0 = 0$. We denote to the extended complex to the left of $X$ in $C_{n+1}(\mathrm{proj}\Lambda)$ by $_EX$, where $_EX^1 = X'^0$, $_EX^i = X^{i-1}$ for all $i \in \{2, \ldots, n+1\}$, $d_{_EX}^1 = d'^0$, $d_{_EX}^i = d_X^{i-1}$ for all $i \in \{2, \ldots, n\}$.*

  ii. *We say that $X$ can be **extended to the right** in $C_{n+1}(\mathrm{proj}\Lambda)$ if there is a projective $\Lambda$-module $X'^{n+1}$ and a non-zero morphism $d'^n \colon X^n \to X'^{n+1}$ such that $d'^n d_X^{n-1} = 0$. We denote to the extended complex to the right of $X$ in $C_{n+1}(\mathrm{proj}\Lambda)$ by $X_E$, where $X_E^i = X^i$ for all $i \in \{1, \ldots, n\}$, $X_E^{n+1} = X'^{n+1}$, $d_{X_E}^i = d_X^i$ for all $i \in \{1, \ldots, n-1\}$, $d_{X_E}^n = d'^n$.*

**Remark 2.6.** Let $X$ be an indecomposable complex in $C_n(\mathrm{proj}\Lambda)$.

  i. We noted that if $X$ can be extended to the left in $C_{n+1}(\mathrm{proj}\Lambda)$ (resp. if $X$ can be extended to the right in $C_{n+1}(\mathrm{proj}\Lambda)$, then we can take $X'^0$ as an indecomposable projective $\Lambda$-module (resp. then we can take $X'^{n+1}$ as an indecomposable projective $\Lambda$-module).

ii. If $X$ can be extended to the left in $C_{n+1}(\text{proj}\Lambda)$ and $X'^0$ is an indecomposable projective $\Lambda$-module, then $_EX$ is an indecomposable complex in $C_{n+1}(\text{proj}\Lambda)$. In fact, there is $X'^0$ an indecomposable projective $\Lambda$-module and non-zero morphism $d'^0\colon X'^0 \to X^1$ such that $d_X^1 d'^0 = 0$.

1. If $X^2 \neq 0$, then we consider the following morphism $f$ in $K^b(\text{proj}\Lambda)$

$$\begin{array}{ccccccccccc}
X'\colon \cdots & \longrightarrow & 0 & \longrightarrow & X'^0 & \longrightarrow & 0 & \longrightarrow & \cdots \longrightarrow & 0 & \longrightarrow & 0 & \longrightarrow \cdots \\
& & \downarrow f & & \downarrow d'^0 & & \downarrow & & & \downarrow & & \downarrow & \\
X\colon \cdots & \longrightarrow & 0 & \longrightarrow & X^1 & \xrightarrow{d_X^1} & X^2 & \longrightarrow & \cdots \longrightarrow & X^n & \longrightarrow & 0 & \longrightarrow \cdots
\end{array}$$

We note that $f$ is non-zero and since $X$ is an indecomposable complex thus $f$ is not invertible. As $\text{Hom}_{K^b(\text{proj}\Lambda)}(X, X'[1]) = 0$ and $X'$ is an indecomposable complex, D. Happel and D. Zacharia in [[19], Corollary 1.4], showed that $C_f$ is an indecomposable complex, where $C_f\colon \cdots \to 0 \to X'^0 \to X^1 \to \cdots \to X^n \to 0 \to \cdots$. Finally, we noted that $C_f[-1] \cong {}_EX$.

2. If $X^2 = 0$, then $X^3 = \cdots = X^n = 0$. Now, if $d'^0$ is an isomorphism, then $_ES(X^1) \cong J_1(X^1)$. And if $d'^0$ is not isomorphism, then an argument similar to 1, is used, so $_ES(X^1) \cong C_f[-1] = \cdots 0 \to X'^0 \to X^1 \to 0 \to \cdots$ is an indecomposable complex.

iii. If $X$ can be extended to the right in $C_{n+1}(\text{proj}\Lambda)$ and $X'^{m+1}$ is an indecomposable projective $\Lambda$-module, then $X_E$ is an indecomposable complex in $C_{n+1}(\text{proj}\Lambda)$. This statement is proved similarly to ii.

iv. C. Chaio, A.G. Chaio, and I. Pratti in [[10], Lemma 4.2], show that $X$ can be extended to the left in $C_{n+1}(\text{proj}\Lambda)$ if and only if $d_X^1$ is not a monomorphism and $X$ can be extended to the right in $C_{n+1}(\text{proj}\Lambda)$ if and only if $\text{Hom}_\Lambda(\text{coker} d_X^{n-1}, \Lambda) \neq 0$.

Also in [10], C. Chaio, A.G. Chaio, and I. Pratti, showed the following result, which has been modified, using the functors defined before Remark 2.3

**Proposition 2.7.** *[[10], Lemma 4.3] Let $X$ and $Y$ be complexes in $C_n(\text{proj}\Lambda)$. If $f\colon X \to Y$ be an irreducible morphism in $C_n(\text{proj}\Lambda)$ with $n \geq 2$, then we have the following.*

i. *If $X$ is an indecomposable and it can not be extended to the left in $C_{n+1}(\text{proj}\Lambda)$, then $i_{-1}^n(f)\colon i_{-1}^n(X) \to i_{-1}^n(Y)$ is an irreducible morphism in $C_{n+1}(\text{proj}\Lambda)$.*

ii. *If $Y$ is an indecomposable and it can not be extended to the right in $C_{n+1}(\text{proj}\Lambda)$, then $i_0^n(f)\colon i_0^n(X) \to i_0^n(Y)$ is an irreducible morphism in $C_{n+1}(\text{proj}\Lambda)$.*

The following proposition, which is the result of joining Corollary 2 and Proposition 3 was proved by H. Giraldo and H. Merklen in [16].

**Proposition 2.8.** *Let $X$ and $Y$ be complexes in $C_n(\text{proj}\Lambda)$. If $f\colon X \to Y$ be an irreducible morphism in $C_n(\text{proj}\Lambda)$, one of the next conditions holds:*

i. *$f^j$ is a section for all $j \in \{1, \ldots, n\}$.*

ii. *$f^j$ is a retraction for all $j \in \{1, \ldots, n\}$.*

iii. there is an $i_0 \in \{1, \ldots, n\}$ such that $f^{i_0}$ is irreducible, the morphisms $f^j$ are sections for all $j > i_0$ and the morphisms $f^j$ are retractions for all $j < i_0$ with $j \in \{1, \ldots, n\}$.

Now, we will provide some basic results and remarks for indecomposable complexes in $\mathrm{K}^b(\mathrm{proj}\Lambda)$.

**Lemma 2.9.** *Let $X$ and $Y$ be indecomposable complexes in $\mathrm{C}_n(\mathrm{proj}\Lambda)$ where $X$ and $Y$ are not $\mathrm{J}_k(P)$ for all $k \in \{1, \ldots, n-1\}$ form being $P$ some indecomposable projective $\Lambda$-module. If $f: X \to Y$ is an irreducible morphism in $\mathrm{C}_n(\mathrm{proj}\Lambda)$ with $X^1, Y^n \neq 0$, then $\mathrm{C}_f$ is an indecomposable in $\mathrm{K}^b(\mathrm{proj}\Lambda)$ and $\ell(\mathrm{C}_f) = n$.*

**Proof.** Using Theorem 5 stated by H. Giraldo and H. Merklen in [16], we have that $X$ and $Y$ are minimal projective complexes and by Theorem 6 stated in [16], $f$ is an irreducible morphism in $\mathrm{K}^b(\mathrm{proj}\Lambda)$; therefore, by the first proposition stated by E.R. Alvares, S.M. Fernandes, and H. Giraldo in [1], section 5.1, the cone of $f$ is indecomposable in $\mathrm{K}^b(\mathrm{proj}\Lambda)$, since $X^1 \neq 0$ and $Y^n \neq 0$, then $\ell(\mathrm{C}_f) = n$. □

From here on, $\Lambda$ will denote a piecewise hereditary algebra or $\Lambda$ will denote an algebra that has finite strong global dimension (by Theorem 2.2).

**Remark 2.10.**

i. If $X \in \mathrm{K}^b(\mathrm{proj}\Lambda)$ is an indecomposable complex with $\ell(X) = \eta$, then X is not (up to traslation) either indecomposable $\mathcal{E}_{\eta+1}$-projective or indecomposable $\mathcal{E}_{\eta+1}$-injective.

   In fact, suppose $X$ is either indecomposable $\mathcal{E}_{\eta+1}$-projective or indecomposable $\mathcal{E}_{\eta+1}$-injective, then by Proposition 2.1 $X$ is some of the following complexes in $\mathrm{C}_{\eta+1}(\mathrm{proj}\Lambda)$:

   (1)  $\mathrm{T}(P)$,   (2)  $\mathrm{S}(P)$,   or   (3)  $\mathrm{J}_k(P)$ con $k \in \{1, \ldots, \eta\}$,

   being $P$ an indecomposable projective $\Lambda$-module. To the complex $\mathrm{T}(P)$ and $\mathrm{S}(P)$, we have that $\ell(\mathrm{T}(P)) = 0 = \ell(\mathrm{S}(P))$, and as complex $\mathrm{J}_k(P)$ with $k \in \{1, \ldots, \eta\}$ is zero in $\mathrm{K}^b(\mathrm{proj}\Lambda)$, then $\ell(\mathrm{J}_k(P)) = 0$. Therefore, $\ell(X) = 0$.

ii. If $X$ is a non-zero indecomposable complex in $\mathrm{K}^b(\mathrm{proj}\Lambda)$ with $\ell(X) = \eta$, then there are $r, s \in \mathbb{Z}$ with $r < s$ such that $X^r \neq 0$ and $X^s \neq 0$, $X^i = 0$ for $i > s$ and $i < r$, and $X^j \neq 0$ for $j \in \{r, r+1, \ldots, s\}$. So, $X[p] \in \mathrm{C}_{\eta+1}(\mathrm{proj}\Lambda)$ for some $p \in \mathbb{Z}$.

iii. If $X$ is an indecomposable complex in $\mathrm{K}^b(\mathrm{proj}\Lambda)$ with $\ell(X) = \eta$, then the following occurs: $X[p] \in \mathrm{C}_{\eta+i}(\mathrm{proj}\Lambda)$ with $i \geq 1$ for some $p \in \mathbb{Z}$, and $X[p] \notin \mathrm{C}_\eta(\mathrm{proj}\Lambda)$.

## 3. Main results

In the first part of this section, we will consider some results about the shape of the indecomposable complexes and almost split $\mathcal{E}_{\eta+i}$-sequences in $\mathrm{C}_{\eta+i}(\mathrm{proj}\Lambda)$ with $i \geq 2$. In the second part, we will prove that, up to translation, all irreducible morphisms and all almost split $\mathcal{E}_{\eta+1}$-sequences in $\mathrm{C}_{\eta+1}(\mathrm{proj}\Lambda)$ are irreducible morphisms and almost split $\mathcal{E}_{\eta+i}$-sequences in $\mathrm{C}_{\eta+i}(\mathrm{proj}\Lambda)$ for $i \geq 2$ respectively (see Theorem 3.21). Finally, Theorem 3.22 and Algorithm 3.24, allow us to find the finite strong global dimension.

**Lemma 3.1.** *Every indecomposable complex in* $C_{\eta+2}(\text{proj}\Lambda)$ *is a complex such that either the cell one or cell $\eta+2$ is zero.*

**Proof.** Let $X$ be an indecomposable complex in $C_{\eta+2}(\text{proj}\Lambda)$, then $X = J_k(P)$ with $P$ being an indecomposable projecive $\Lambda$-module for $k \in \{1,\ldots,\eta+1\}$ or $X \neq J_k(P)$.

If $X = J_k(P)$, then $X$ satisfies the condition.

If $X \neq J_k(P)$, then X does not have summands of the form $J_k(P)$. Therefore, by H. Giraldo and H. Merklen in [[16], Theorem 5] (also see [[6], Proposition 6.6), $X$ is also indecomposable in $K^b(\text{proj}\Lambda)$. If we suppose that $X^1 \neq 0$ and $X^{\eta+2} \neq 0$, then $\ell(X) = \eta+1$, which contradicts the fact that $\eta = \text{s.gl.dim}(\Lambda)$. □

**Corollary 3.2.** *Every indecomposable complex in* $C_{\eta+i}(\text{proj}\Lambda)$ *is a complex such that either the cell one or cell $\eta+i$ is zero for $i \geq 2$.*

**Proposition 3.3.** *Let $X$ and $Y$ be indecomposable complexes in $C_{\eta+i}(\text{proj}\Lambda)$ with $i \geq 2$. If $f: X \to Y$ is an irreducible morphism in $C_{\eta+i}(\text{proj}\Lambda)$, then $X^1 = 0 = Y^1$ or $X^{\eta+i} = 0 = Y^{\eta+i}$.*

**Proof.** Let $f: X \to Y$ be an irreducible morphism in $C_{\eta+i}(\text{proj}\Lambda)$. Let's suppose that $X^1 \neq 0$ or $Y^1 \neq 0$, then we will show that $X^{\eta+i} = 0 = Y^{\eta+i}$.

$X^1 \neq 0$. If $X^{\eta+i} \neq 0$, then $X^1 = 0$ because $\ell(X) < \eta$; therefore, $X^{\eta+i} = 0$. On the other hand, if $Y^{\eta+i} \neq 0$ then $X$ and $Y$ are not the form $J_k(P)$ for all $k \in \{1,\ldots,\eta+i-1\}$ with $P$ being some indecomposable projective $\Lambda$-module, since in other way $f = 0$. By Lemma 2.9, we obtained that $\ell(C_f) \geq \eta+2$, which contradicts the fact that $\eta = \text{s.gl.dim}(\Lambda)$, thus $Y^{\eta+i} = 0$.

$Y^1 \neq 0$. If $Y^{\eta+i} \neq 0$, then $Y^1 = 0$ because $\ell(Y) < \eta$; therefore, $Y^{\eta+i} = 0$. Now, if $X^{\eta+i} \neq 0$, the morphism $f^{\eta+i}: X^{\eta+i} \to 0$ is not a section in proj$\Lambda$ and the morphism $f^1: 0 \to Y^1$ is not a retraction in proj$\Lambda$, which is a contradiction by Proposition 2.8. □

**Corollary 3.4.** *Let $X$ and $Y$ be indecomposable complexes in $C_{\eta+i}(\text{proj}\Lambda)$ with $i \geq 2$ and $f: X \to Y$ an irreducible morphism in $C_{\eta+i}(\text{proj}\Lambda)$. The following conditions hold:*

  i. *If $\ell(X) = l - j = \eta$, then $Y^{j-r} = Y^{l+r} = 0$ for all $r \geq 1$.*
  ii. *If $\ell(Y) = l - j = \eta$, then $X^{j-r} = X^{l+r} = 0$ for all $r \geq 1$.*

That is, with $\eta = l - j$, so if $\ell(X) = \eta$ or $\ell(Y) = \eta$ then $f$ has the following form:

$$\begin{array}{ccccccccccccc}
X: & 0 & \longrightarrow & \cdots & \longrightarrow & 0 & \longrightarrow & X^j & \longrightarrow & \cdots & \longrightarrow & X^l & \longrightarrow 0 \longrightarrow \cdots \longrightarrow 0 \to 0 \to \cdots \\
& & & & & & & \downarrow f^j & & & & \downarrow f^l & \\
Y: & 0 & \longrightarrow & \cdots & \longrightarrow & 0 & \longrightarrow & Y^j & \longrightarrow & \cdots & \longrightarrow & Y^l & \longrightarrow 0 \longrightarrow \cdots \longrightarrow 0 \to 0 \to \cdots.
\end{array}$$

with $\downarrow f$ on the left.

**Proof.** We obtained a) and b) in a way similar to what was carried out in the proof of Proposition 3.3. □

**Corollary 3.5.** *Let $X$ and $Y$ be indecomposable complexes in $C_{\eta+2}(\text{proj}\Lambda)$. If $f: X \to Y$ is an irreducible morphism in $C_{\eta+2}(\text{proj}\Lambda)$, then there is an irreducible morphism $g$ in $C_{\eta+1}(\text{proj}\Lambda)$ such that $i^{\eta+1}_{-1}(g) = f$ or $i^{\eta+1}_0(g) = f$.*

**Proof.** By Proposition 3.3 $X^1 = 0 = Y^1$ or $X^{\eta+2} = 0 = Y^{\eta+2}$.

If $X^1 = 0 = Y^1$, then we consider $g = \pi_1^{\eta+1}(f)$ and by Lemma 2.4, $g$ is an irreducible morphism in $C_{\eta+1}(\text{proj}\Lambda)$ and by Remark 2.3 ii. $i_{-1}^{\eta+1}(g) = f$.

If $X^{\eta+2} = 0 = Y^{\eta+2}$, then we consider $g = \pi_0^{\eta+1}(f)$ by Lemma 2.4, $g$ is an irreducible morphism in $C_{\eta+1}(\text{proj}\Lambda)$ and by Remark 2.3 iii. $i_0^{\eta+1}(g) = f$. $\square$

**Remark 3.6.** From Corollary 3.5, if $X$ and $Y$ are indecomposable complexes in $C_{\eta+2}(\text{proj}\Lambda)$ and $f: X \to Y$ is an irreducible morphism in $C_{\eta+2}(\text{proj}\Lambda)$, then up to functor $i_{-1}^{\eta+1}$ or $i_0^{\eta+1}$ - which are translation and isomorphisms by Remark 2.3 - we can say that, up to translation, $f$ is an irreducible morphism in $C_{\eta+1}(\text{proj}\Lambda)$.

**Lemma 3.7.** If $\delta: X \xrightarrow{f} Y \xrightarrow{g} Z$ is an almost split $\mathcal{E}_{\eta+2}$-sequence in $C_{\eta+2}(\text{proj}\Lambda)$; then we have the following:

i. If $X^1 = 0 = Y^1$ and $Z^1 = 0$, then $\pi_1^{\eta+1}(\delta): \pi_1^{\eta+1}(X) \to \pi_1^{\eta+1}(Y) \to \pi_1^{\eta+1}(Z)$ is an almost split $\mathcal{E}_{\eta+1}$-sequence in $C_{\eta+1}(\text{proj}\Lambda)$.

ii. If $X^{\eta+2} = 0 = Y^{\eta+2}$ and $Z^{\eta+2} = 0$, then $\pi_0^{\eta+1}: \pi_0^{\eta+1}(X) \to \pi_0^{\eta+1}(Y) \to \pi_0^{\eta+1}(Z)$ is an almost split $\mathcal{E}_{\eta+1}$-sequence in $C_{\eta+1}(\text{proj}\Lambda)$.

**Proof.** Let $\delta: X \xrightarrow{f} Y \xrightarrow{g} Z$ be an almost split sequence in $C_{\eta+2}(\text{proj}\Lambda)$.

i. Let's suppose that $X^1 = 0 = Y^1$ and $Z^1 = 0$, since $f$ is a minimal left almost split morphism in $C_{\eta+2}(\text{proj}\Lambda)$, then $f$ is an irreducible morphism in $C_{\eta+2}(\text{proj}\Lambda)$. By the Lemma 2.4, $\pi_1^{\eta+1}(f): \pi_1^{\eta+1}(X) \to \pi_1^{\eta+1}(Y)$ is an irreducible morphism in $C_{\eta+1}(\text{proj}\Lambda)$.

We will prove that $\pi_1^{\eta+1}(f): \pi_1^{\eta+1}(X) \to \pi_1^{\eta+1}(Y)$ and $\pi_1^{\eta+1}(g): \pi_1^{\eta+1}(Y) \to \pi_1^{\eta+1}(Z)$ are a minimal left and right almost split morphisms in $C_{\eta+1}(\text{proj}\Lambda)$ respectively.

1. Since $\pi_1^{\eta+1}(f): \pi_1^{\eta+1}(X) \to \pi_1^{\eta+1}(Y)$ is an irreducible morphism in $C_{\eta+1}(\text{proj}\Lambda)$, then $\pi_1^{\eta+1}(f)$ is not a section in $C_{\eta+1}(\text{proj}\Lambda)$.

2. Let $u: \pi_1^{\eta+1}(X) \to U$ be a morphism in $C_{\eta+1}(\text{proj}\Lambda)$ such that $u$ is not a section in $C_{\eta+1}(\text{proj}\Lambda)$, we will prove that there is $v: \pi_1^{\eta+1}(Y) \to U$ a morphism in $C_{\eta+1}(\text{proj}\Lambda)$ such that $v\pi_1^{\eta+1}(f) = u$. We considered the complex $i_{-1}^{\eta+1}(U)$ and the morphism $i_{-1}^{\eta+1}(u)$ in $C_{\eta+2}(\text{proj}\Lambda)$ and we realized that $i_{-1}^{\eta+1}(u)$ is not a section in $C_{\eta+2}(\text{proj}\Lambda)$, since $u$ is not a section in $C_{\eta+1}(\text{proj}\Lambda)$. As $f$ is a minimal left almost split morphism in $C_{\eta+2}(\text{proj}\Lambda)$, there is $h: Y \to i_{-1}^{\eta+1}(U)$, a morphism in $C_{\eta+2}(\text{proj}\Lambda)$, such that $hf = i_{-1}^{\eta+1}(u)$. As $Y^1 = 0 = i_{-1}^{\eta+1}(U)^1$, we have $\pi_1^{\eta+1}(h)\pi_1^{\eta+1}(f) = u$ in $C_{\eta+1}(\text{proj}\Lambda)$.

3. Let $r \in \text{End}(\pi_1^{\eta+1}(Y))$ be such that $r\pi_1^{\eta+1}(f) = \pi_1^{\eta+1}(f)$ in $C_{\eta+1}(\text{proj}\Lambda)$, so we have $i_{-1}^{\eta+1}(r) \in \text{End}(Y)$ and $i_{-1}^{\eta+1}(r)f = f$ in $C_{\eta+2}(\text{proj}\Lambda)$, and since $f$ is minimal in $C_{\eta+2}(\text{proj}\Lambda)$, then $i_{-1}^{\eta+1}(r) \in \text{Aut}(Y)$ in $C_{\eta+2}(\text{proj}\Lambda)$ and, therefore, $r \in \text{Aut}(\pi_1^{\eta+1}(Y))$ in $C_{\eta+1}(\text{proj}\Lambda)$.

Considering 1, 2 and 3, we concluded that, $\pi_1^{\eta+1}(f)$ is a minimal left almost split morphism in $C_{\eta+1}(\text{proj}\Lambda)$. In a similar way, we can prove that $\pi_1^{\eta+1}(g): \pi_1^{\eta+1}(Y) \to \pi_1^{\eta+1}(Z)$ is minimal right almost split morphism in $C_{\eta+1}(\text{proj}\Lambda)$. This proves that $\pi_1^{\eta+1}(\delta): \pi_1^{\eta+1}(X) \to \pi_1^{\eta+1}(Y) \to \pi_1^{\eta+1}(Z)$ is an almost split $\mathcal{E}_{\eta+1}$-sequence in $C_{\eta+1}(\text{proj}\Lambda)$.

ii. The proof is similar to i. □

**Corollary 3.8.** *Let $\delta: X \to Y \to Z$ be an almost split $\mathcal{E}_{\eta+i}$-sequence in $C_{\eta+i}(\text{proj}\Lambda)$ with $i \geq 2$, then we have the following:*

i. *If $X^1 = 0 = Y^1$ and $Z^1 = 0$, then $\pi_1^{\eta+i}(\delta): \pi_1^{\eta+i}(X) \to \pi_1^{\eta+i}(Y) \to \pi_1^{\eta+i}(Z)$ is an almost split $\mathcal{E}_{\eta+(i-1)}$-sequence in $C_{\eta+(i-1)}(\text{proj}\Lambda)$.*

ii. *If $X^{\eta+i} = 0 = Y^{\eta+i}$ and $Z^{\eta+i} = 0$, then $\pi_0^{\eta+i}: \pi_0^{\eta+i}(X) \to \pi_0^{\eta+i}(Y) \to \pi_0^{\eta+i}(Z)$ is an almost split $\mathcal{E}_{\eta+(i-1)}$-sequence in $C_{\eta+(i-1)}(\text{proj}\Lambda)$.*

**Proposition 3.9.** *If $\delta: X \xrightarrow{f} Y \xrightarrow{g} Z$ is an almost split $\mathcal{E}_{\eta+2}$-sequence in $C_{\eta+2}(\text{proj}\Lambda)$, then there is an almost split $\mathcal{E}_{\eta+1}$-sequences $\delta'$ in $C_{\eta+1}(\text{proj}\Lambda)$ such that $i_{-1}^{\eta+1}(\delta') = \delta$ or $i_0^{\eta+1}(\delta') = \delta$.*

**Proof.** Since $X \xrightarrow{f} Y$ is an irreducible morphism in $C_{\eta+2}(\text{proj}\Lambda)$ by Proposition 3.3 we have that $X^1 = 0 = Y^1$ or $X^{\eta+2} = 0 = Y^{\eta+2}$.

If $X^1 = 0 = Y^1$, then $Z^1 = 0$ since $0 \to X^1 \to Y^1 \to Z^1 \to 0$ is split. Let's consider $\delta' = \pi_1^{\eta+1}(\delta)$, by Lemma 3.7 $\pi_1^{\eta+1}(\delta)$ is an almost split $\mathcal{E}_{\eta+1}$-sequence in $C_{\eta+1}(\text{proj}\Lambda)$ and by Remark 2.3 ii. $i_{-1}^{\eta+1}(\delta') = \delta$.

If $X^{\eta+2} = 0 = Y^{\eta+2}$, then $Z^{\eta+2} = 0$. Let's consider $\delta' = \pi_0^{\eta+1}(\delta)$, by Lemma 3.7 $\pi_0^{\eta+1}(\delta)$ is an almost split $\mathcal{E}_{\eta+1}$-sequence in $C_{\eta+1}(\text{proj}\Lambda)$ and by Remark 2.3 iii. $i_0^{\eta+1}(\delta') = \delta$. □

**Remark 3.10.** According to Proposition 3.9, if $\delta: X \xrightarrow{f} Y \xrightarrow{g} Z$ is an almost split $\mathcal{E}_{\eta+1}$-sequence in $C_{\eta+2}(\text{proj}\Lambda)$, then up to functor $i_{-1}^{\eta+1}$ or $i_0^{\eta+1}$ - which are translation and isomorphisms by Remark 2.3 - we can say that up to translation, $\delta$ is an almost split $\mathcal{E}_{\eta+1}$-sequence in $C_{\eta+1}(\text{proj}\Lambda)$.

**Lemma 3.11.** *Let $X$ and $Y$ be indecomposable complexes in $C_{\eta+1}(\text{proj}\Lambda)$. If $f: X \to Y$ is an irreducible morphism in $C_{\eta+1}(\text{proj}\Lambda)$, then we have the following:*

i. *If $X$ can be extended to the left in $C_{\eta+2}(\text{proj}\Lambda)$, then $Y$ can not be extended to the right in $C_{\eta+2}(\text{proj}\Lambda)$.*

ii. *If $Y$ can be extended to the right in $C_{\eta+2}(\text{proj}\Lambda)$, then $X$ can not be extended to the left in $C_{\eta+2}(\text{proj}\Lambda)$.*

**Proof.**

i. If $X$ can be extended to the left in $C_{\eta+2}(\text{proj}\Lambda)$, then $X^1 \neq 0$. Now, we suppose that $Y$ can be extended to the right in $C_{\eta+2}(\text{proj}\Lambda)$; therefore, $Y^{\eta+1} \neq 0$ and also $X$ and $Y$ are not the form $J_k(P)$ for all $k \in \{1, \ldots, \eta\}$ with $P$ being some indecomposable projective $\Lambda$-module. By Lemma 2.9, we have that $\ell(C_f) = \eta + 1$, which contradicts the fact that $\eta = \text{s.gl.dim}(\Lambda)$.

ii. The proof is similar to i. □

**Proposition 3.12.** *Let $X$ and $Y$ be indecomposable complexes in $C_{\eta+1}(\mathrm{proj}\Lambda)$. If $f\colon X\to Y$ is an irreducible morphism in $C_{\eta+1}(\mathrm{proj}\Lambda)$, then there is an irreducible morphism $g$ in $C_{\eta+2}(\mathrm{proj}\Lambda)$ such that $\pi_1^{\eta+1}(g)=f$ or $\pi_0^{\eta+1}(g)=f$.*

**Proof.**

If $X$ can not be extended to the left in $C_{\eta+2}(\mathrm{proj}\Lambda)$, then let's consider $g=i_{-1}^{\eta+1}(f)$ by Proposition 2.7, $i_{-1}^{\eta+1}(f)$ is an irreducible morphism in $C_{\eta+2}(\mathrm{proj}\Lambda)$ and by Remark 2.3 ii. $\pi_1^{\eta+1}(g)=f$.

If $X$ can be extended to the left in $C_{\eta+2}(\mathrm{proj}\Lambda)$, then $Y$ can not be extended to the right in $C_{\eta+2}(\mathrm{proj}\Lambda)$ by Lemma 3.11, then let's consider $g=i_0^{\eta+1}(f)$, by Proposition 2.7, $i_0^{\eta+1}(f)$ is an irreducible morphism in $C_{\eta+2}(\mathrm{proj}\Lambda)$ and by Remark 2.3 iii. $\pi_0^{\eta+1}(g)=f$. □

**Remark 3.13.** According to Proposition 3.12, if $X$ and $Y$ are indecomposable complexes in $C_{\eta+1}(\mathrm{proj}\Lambda)$ and $f\colon X\to Y$ is an irreducible morphism in $C_{\eta+1}(\mathrm{proj}\Lambda)$, then up to functor $\pi_1^{\eta+1}$ or $\pi_0^{\eta+1}$ - which are translation and isomorphisms by Remark 2.3 - we can say that up to translation, $f$ is an irreducible morphism in $C_{\eta+2}(\mathrm{proj}\Lambda)$.

**Lemma 3.14.** *Let $X\in C_{\eta+1}(\mathrm{proj}\Lambda)$ be an indecomposable complex.*

  i. *If $X$ can be extended to the left in $C_{\eta+2}(\mathrm{proj}\Lambda)$, then $X^{\eta+1}=0$.*

  ii. *If $X$ can be extended to the right in $C_{\eta+2}(\mathrm{proj}\Lambda)$, then $X^1=0$.*

**Proof.** Let $X\in C_{\eta+1}(\mathrm{proj}\Lambda)$ be an indecomposable complex. For both cases i and ii, we noted that $X\neq J_k(P)$ with $k\in\{1,\ldots,\eta\}$ and $P$ is an indecomposable projective $\Lambda$-module.

  i. We assumed that $X^{\eta+1}\neq 0$ and $X$ can be extended to the left in $C_{\eta+2}(\mathrm{proj}\Lambda)$, so there is $X'^0\in\mathrm{proj}\Lambda$ and non-zero morphism $d'^0\colon X'^0\to X^1$ such that $d_X^1 d'^0=0$. We considered the following morphism $f$ in $K^b(\mathrm{proj}\Lambda)$:

$$\begin{array}{ccccccccccc}
X'\colon & \cdots\longrightarrow & 0 & \longrightarrow & X'^0 & \longrightarrow & 0 & \longrightarrow\cdots\longrightarrow & 0 & \longrightarrow & 0 \longrightarrow 0 \longrightarrow\cdots \\
& & \downarrow f & & \downarrow d'^0 & & \downarrow & & \downarrow & & \downarrow \quad\downarrow \\
X\colon & \cdots\longrightarrow & 0 & \longrightarrow & X^1 & \xrightarrow{d_X^1} & X^2 & \longrightarrow\cdots\longrightarrow & X^\eta & \xrightarrow{d_X^\eta} & X^{\eta+1} \longrightarrow 0 \longrightarrow\cdots
\end{array}$$

  We noted that $f$ is non-zero and not invertible, since $X$ is indecomposable. As $\mathrm{Hom}_{K(\mathrm{proj}\Lambda)}(X,X'[1])=0$ according to D. Happel and D. Zacharia in [[19], Corollary 1.4], then $C_f$ is indecomposable and $\ell(C_f)=\eta+1$. But this is a contradiction, because $\eta=\mathrm{s.gl.dim}(\Lambda)$.

  ii. Similar to i. □

**Lemma 3.15.** *Let $\delta\colon X\xrightarrow{f}Y\xrightarrow{g}Z$ be an almost split $\mathcal{E}_{\eta+1}$-sequence in $\mathrm{C}_{\eta+1}(\mathrm{proj}\Lambda)$, then we have the following:*

i. *If $X$ can be extended to the left in $\mathrm{C}_{\eta+2}(\mathrm{proj}\Lambda)$, then $Z$ can not be extended to the right in $\mathrm{C}_{\eta+2}(\mathrm{proj}\Lambda)$.*

ii. *If $Z$ can be extended to the right in $\mathrm{C}_{\eta+2}(\mathrm{proj}\Lambda)$, then $X$ can not be extended to the left in $\mathrm{C}_{\eta+2}(\mathrm{proj}\Lambda)$.*

**Proof.**

i. Suppose that $Z$ can be extended to the right in $\mathrm{C}_{\eta+2}(\mathrm{proj}\Lambda)$, then $Z^{\eta+1}\neq 0$ and $Z\neq \mathrm{J}_k(P)$ for all $k\in\{2,\ldots,\eta\}$ with $P$ some indecomposable projective $\Lambda$-module; furthermore, by Lemma 3.14, $Z^1=0$. Since $0\to X^1 \to Y^1 \to Z^1 \to 0$ is split, then $X^1\cong Y^1$. $X$ can be extended to the left in $\mathrm{C}_{\eta+2}(\mathrm{proj}\Lambda)$ so $X^1\neq 0$ and, therefore, $Y^1\neq 0$, thus there is $k\in\{1,...,m\}$ such that $Y_k^1\neq 0$, where $Y=\bigoplus_{i=1}^{m}Y_i$ with $Y_i$ being an indecomposable complex in $\mathrm{C}_{\eta+1}(\mathrm{proj}\Lambda)$ for all $i\in\{1,\ldots,m\}$. Since $X\xrightarrow{f}Y\xrightarrow{g}Z$ is an almost split $\mathcal{E}_{\eta+1}$-sequence in $\mathrm{C}_{\eta+1}(\mathrm{proj}\Lambda)$ by R. Bautista in [[4], Remark 2.1], the morphism $g=[\,g_1\ \cdots\ g_k\ \cdots\ g_m\,]$ is an irreducible morphism in $\mathrm{C}_{\eta+1}(\mathrm{proj}\Lambda)$ and also a radical morphism in $\mathrm{C}_{\eta+1}(\mathrm{proj}\Lambda)$ from [[5], Lemma 2.16]. Therefore, by R. Bautista and M.J. Souto Salorio in [[5], Proposition 2.18], we have that $g_i$ is an irreducible morphism in $\mathrm{C}_{\eta+1}(\mathrm{proj}\Lambda)$ for all $i\in\{1,\ldots,m\}$. We considered $g_k\colon Y_k\to Z$ in $\mathrm{C}_{\eta+1}(\mathrm{proj}\Lambda)$, so

$$\begin{array}{ccccccccccc}
Y_k\colon & \cdots\longrightarrow & 0 & \longrightarrow & Y_k^1 & \xrightarrow{d_Y^1} & Y_k^2 & \longrightarrow & \cdots & \longrightarrow & Y_k^\eta & \xrightarrow{d_Y^\eta} & Y_k^{\eta+1} & \longrightarrow 0 \longrightarrow \cdots\\
 & & \downarrow & & \downarrow & & \downarrow 0 & & \downarrow g_k^2 & & & & \downarrow & & \downarrow & & \downarrow \\
Z\colon & \cdots\longrightarrow & 0 & \longrightarrow & 0 & \xrightarrow{0} & Z^2 & \longrightarrow & \cdots & \longrightarrow & Z^\eta & \xrightarrow{d_Z^\eta} & Z^{\eta+1} & \longrightarrow 0 \longrightarrow \cdots
\end{array}$$

1. Let $Y_k\neq \mathrm{J}_1(Y_k^1)$, since $Y_k^1\neq 0$ and $Z^{\eta+1}\neq 0$. By Lemma 2.9, we obtained that $\ell(\mathrm{C}_{g_k})=\eta+1$, contradicting the fact that $\eta=\mathrm{s.gl.dim}(\Lambda)$.

2. Let $Y_k=\mathrm{J}_1(Y_k^1)$. In this case $g_k^2=0$. By Proposition 2.8, $g_k^2\colon Y_k^1\to Z^2$ is an irreducible morphism in $\mathrm{proj}\Lambda$; therefore, $g_k^2$ is a monomorphism since $Z^2$ is a projective $\Lambda$-module. So $0=\ker g_k^2=Y_k^1$, contradicting the fact that $Y_k^1\neq 0$.

All of the above allows us to conclude that, $Z$ can not be extended to the right in $\mathrm{C}_{\eta+2}(\mathrm{proj}\Lambda)$.

ii. The proof is similar to i. $\square$

**Lemma 3.16.** *If $\delta\colon X\xrightarrow{f}Y\xrightarrow{g}Z$ is an almost split $\mathcal{E}_{\eta+1}$-sequence in $\mathrm{C}_{\eta+1}(\mathrm{proj}\Lambda)$, then we have the following:*

i. *If $X$ can not be extended to the left in $\mathrm{C}_{\eta+2}(\mathrm{proj}\Lambda)$, then $i_{-1}^{\eta+1}(\delta)$ is an almost split $\mathcal{E}_{\eta+2}$-sequence in $\mathrm{C}_{\eta+2}(\mathrm{proj}\Lambda)$.*

ii. *If $X$ can be extended to the left $\mathrm{C}_{\eta+2}(\mathrm{proj}\Lambda)$, then $i_{0}^{\eta+1}(\delta)$ is an almost split $\mathcal{E}_{\eta+2}$-sequence in $\mathrm{C}_{\eta+2}(\mathrm{proj}\Lambda)$.*

**Proof.** Let $\delta\colon X \xrightarrow{f} Y \xrightarrow{g} Z$ be an almost split $\mathcal{E}_{\eta+1}$-sequence in $C_{\eta+1}(\operatorname{proj}\Lambda)$.

  i. Suppose that $X$ can not be extended to the left in $C_{\eta+2}(\operatorname{proj}\Lambda)$. As $f$ is an irreducible morphism in $C_{\eta+1}(\operatorname{proj}\Lambda)$, by Proposition 2.7, $i_{-1}^{\eta+1}(f)\colon i_{-1}^{\eta+1}(X) \to i_{-1}^{\eta+1}(Y)$ is an irreducible morphism in $C_{\eta+2}(\operatorname{proj}\Lambda)$. We will prove that $i_{-1}^{\eta+1}(\delta)$ is an almost split $\mathcal{E}_{\eta+2}$-sequence in $C_{\eta+2}(\operatorname{proj}\Lambda)$; for this, we will show that $i_{-1}^{p}(f)$ is a minimal left almost split morphims in $C_{\eta+2}(\operatorname{proj}\Lambda)$. In fact, let $v\colon i_{-1}^{\eta+1}(X) \to X'$ be morphism in $C_{\eta+2}(\operatorname{proj}\Lambda)$ such that $v$ is not a section. We considered $X''\colon (X'^2 \to X'^3 \to \cdots \to X'^{\eta+2})$ a complex in $C_{\eta+1}(\operatorname{proj}\Lambda)$ and $v'\colon X \to X''$ in $C_{\eta+1}(\operatorname{proj}\Lambda)$ given by $v'\colon (v^2,\ldots,v^{\eta+2})$. We noted that $v'$ is not a section in $C_{\eta+1}(\operatorname{proj}\Lambda)$, since $X$ can not be extended to the left in $C_{\eta+2}(\operatorname{proj}\Lambda)$ and $v$ is not a section in $C_{\eta+2}(\operatorname{proj}\Lambda)$. As $f$ is a minimal left almost split morphism in $C_{\eta+1}(\operatorname{proj}\Lambda)$, there is $v''\colon Y \to X''$ such that $v''f = v'$. If we define $h\colon i_{-1}^{\eta+1}(Y) \to X'$ by $h\colon (0, v''^1,\ldots,v''^{\eta+1})$, then we have $h i_{-1}^{\eta+1}(f) = v$. Since $f$ is left minimal morphism $C_{\eta+1}(\operatorname{proj}\Lambda)$, we concluded that $i_{-1}^{\eta+1}(f)$ is a minimal left morphism in $C_{\eta+2}(\operatorname{proj}\Lambda)$. Finally, we have that $i_{-1}^{\eta+1}(f)$ is a minimal left almost split morphism in $C_{\eta+2}(\operatorname{proj}\Lambda)$. In a similar way, we can prove that $i_{-1}^{\eta+1}(g)$ is a minimal right almost split morphism in $C_{\eta+2}(\operatorname{proj}\Lambda)$. This proves that $i_{-1}^{\eta+1}(\delta)$ is an almost split $\mathcal{E}_{\eta+2}$-sequence in $C_{\eta+2}(\operatorname{proj}\Lambda)$.

  ii. Suppose that $X$ can be extended to the left in $C_{\eta+2}(\operatorname{proj}\Lambda)$, by Lemma 3.15 $Z$ can not be extended to the right in $C_{\eta+2}(\operatorname{proj}\Lambda)$. Similarly to item i., we proved that $i_{0}^{\eta+1}(g)$ is minimal right almost split morphism in $C_{\eta+2}(\operatorname{proj}\Lambda)$ and $i_{0}^{\eta+1}(\delta)$ is an almost split $\mathcal{E}_{\eta+2}$-sequence in $C_{\eta+2}(\operatorname{proj}\Lambda)$. □

**Proposition 3.17.** *If $\delta\colon X \to Y \to Z$ is an almost split $\mathcal{E}_{\eta+1}$-sequence in $C_{\eta+1}(\operatorname{proj}\Lambda)$, then there is an almost split $\mathcal{E}_{\eta+2}$-sequence $\delta'$ in $C_{\eta+2}(\operatorname{proj}\Lambda)$ such that $\pi_1^{\eta+1}(\delta') = \delta$ or $\pi_0^{\eta+1}(\delta') = \delta$.*

**Proof.** Let $\delta\colon X \xrightarrow{f} Y \xrightarrow{g} Z$ be an almost split $\mathcal{E}_{\eta+1}$-sequence in $C_{\eta+1}(\operatorname{proj}\Lambda)$.

$X$ can not be extended to the left in $C_{\eta+2}(\operatorname{proj}\Lambda)$. We considered $\delta' = i_{-1}^{\eta+1}(\delta)$, by Lemma 3.16, $i_{-1}^{\eta+1}(\delta)$ is an almost split $\mathcal{E}_{\eta+2}$-sequence in $C_{\eta+2}(\operatorname{proj}\Lambda)$ and by Remark 2.3 ii., $\pi_1^{\eta+1}(\delta') = \delta$.

$X$ can be extended to the left in $C_{\eta+2}(\operatorname{proj}\Lambda)$. We considered $\delta' = i_0^{\eta+1}(\delta)$, by Lemma 3.16, $i_0^{\eta+1}(\delta)$ is an almost split $\mathcal{E}_{\eta+2}$-sequence in $C_{\eta+2}(\operatorname{proj}\Lambda)$ and by Remark 2.3 iii., $\pi_0^{\eta+1}(\delta') = \delta$. □

**Remark 3.18.** From Proposition 3.17, if $\delta\colon X \to Y \to Z$ is an almost split $\mathcal{E}_{\eta+1}$-sequence in $C_{\eta+1}(\operatorname{proj}\Lambda)$, then up to functor $\pi_1^{\eta+1}$ or $\pi_0^{\eta+1}$ - which are translation and isomorphisms by Remark 2.3 - we can say that, up to translation, $\delta$ is an almost split $\mathcal{E}_{\eta+2}$-sequence in $C_{\eta+2}(\operatorname{proj}\Lambda)$.

According to Remark 3.6, Remark 3.10, Remark 3.13, and Remark 3.18, we can state, up to translation, the following proposition:

**Proposition 3.19.** *The following conditions hold:*

  i. *Let $X$ and $Y$ be indecomposable complexes in $C_{\eta+2}(\operatorname{proj}\Lambda)$. Up to translation, $f\colon X \to Y$ is an irreducible morphism in $C_{\eta+2}(\operatorname{proj}\Lambda)$ if and only if $f\colon X \to Y$ is an irreducible morphism in $C_{\eta+1}(\operatorname{proj}\Lambda)$.*

  ii. Up to translation, $\delta\colon X \to Y \to Z$ is an almost split $\mathcal{E}_{\eta+2}$-sequence in $C_{\eta+2}(\operatorname{proj}\Lambda)$ if and only if $\delta\colon X \to Y \to Z$ is an almost split $\mathcal{E}_{\eta+1}$-sequences in $C_{\eta+1}(\operatorname{proj}\Lambda)$.

**Proof.**

  i. ($\Longrightarrow$) Corollary 3.5.
     ($\Longleftarrow$) Proposition 3.12.

  ii. ($\Longrightarrow$) Proposition 3.9.
      ($\Longleftarrow$) Proposition 3.17.  □

**Remark 3.20.**

  i. Let $X$ and $Y$ be indecomposable complexes in $C_{\eta+i}(\operatorname{proj}\Lambda)$ with $i \geq 2$. If $f\colon X \to Y$ is an irreducible morphism in $C_{\eta+i}(\operatorname{proj}\Lambda)$, then $X^1 = 0 = Y^1$ or $X^{\eta+i} = 0 = Y^{\eta+i}$ by Proposition 3.3. Now, by Lemma 2.4, up to translation, $f\colon X \to Y$ is an irreducible morphism in $C_{\eta+i-1}(\operatorname{proj}\Lambda)$. If we continue in this way, we can say that, up to translation, $f\colon X \to Y$ is an irreducible morphism in $C_{\eta+1}(\operatorname{proj}\Lambda)$.

  ii. Let $X$ and $Y$ be indecomposable complexes in $C_{\eta+1}(\operatorname{proj}\Lambda)$. If $f\colon X \to Y$ is an irreducible morphism in $C_{\eta+1}(\operatorname{proj}\Lambda)$, by Remark 3.13 $f$ is an irreducible morphism in $C_{\eta+2}(\operatorname{proj}\Lambda)$. If we continue in this way, we can say that, up to translation, $f\colon X \to Y$ is an irreducible morphism in $C_{\eta+i}(\operatorname{proj}\Lambda)$ with $i \geq 2$.

  iii. If $\delta\colon X \to Y \to Z$ is an almost split $\mathcal{E}_{\eta+i}$-sequence in $C_{\eta+i}(\operatorname{proj}\Lambda)$ with $i \geq 2$, by Corollary 3.8, up to translation, $\delta\colon X \to Y \to Z$ is an almost split $\mathcal{E}_{\eta+(i-1)}$-sequence in $C_{\eta+(i-1)}(\operatorname{proj}\Lambda)$. If we continue in this way, we can say that, up to translation, $\delta\colon X \to Y \to Z$ is an almost split $\mathcal{E}_{\eta+1}$-sequence in $C_{\eta+1}(\operatorname{proj}\Lambda)$.

  iv. If $\delta\colon X \to Y \to Z$ is an almost split $\mathcal{E}_{\eta+1}$-sequence in $C_{\eta+1}(\operatorname{proj}\Lambda)$, by Remark 3.18, we can say that, up to translation, $\delta\colon X \to Y \to Z$ is an almost split $\mathcal{E}_{\eta+2}$-sequence in $C_{\eta+2}(\operatorname{proj}\Lambda)$. If we continue in this way, we can say that, up to translation, $\delta\colon X \to Y \to Z$ is an almost split $\mathcal{E}_{\eta+i}$-sequence in $C_{\eta+i}(\operatorname{proj}\Lambda)$ with $i \geq 2$.

According to Remark 3.20, we can state, up to translation, the following theorem:

**Theorem 3.21.** *Let $\Lambda$ be an algebra with finite strong global dimension where $\eta = \operatorname{s.gl.dim}(\Lambda)$. The following conditions hold:*

  i. *Let $X$ and $Y$ be indecomposable complexes in $C_{\eta+i}(\operatorname{proj}\Lambda)$ with $i \geq 2$. Up to translation, $f\colon X \to Y$ is an irreducible morphism in $C_{\eta+i}(\operatorname{proj}\Lambda)$, if and only if, $f\colon X \to Y$ is an irreducible morphism in $C_{\eta+1}(\operatorname{proj}\Lambda)$.*

  ii. *Up to translation, $\delta\colon X \to Y \to Z$ is an almost split $\mathcal{E}_{\eta+i}$-sequence in $C_{\eta+i}(\operatorname{proj}\Lambda)$ with $i \geq 2$ if and only if $\delta\colon X \to Y \to Z$ is an almost split $\mathcal{E}_{\eta+1}$-sequence in $C_{\eta+1}(\operatorname{proj}\Lambda)$.*

We complete this section with the following Theorem, which allows giving the Algorithm to calculate the strong global dimension.

**Theorem 3.22.** *If $\Lambda$ is a finite dimensional $\Bbbk$-algebra, then $\operatorname{s.gl.dim}(\Lambda) = m_0 - 2$ where*

$$m_0 := \min\{n \geq 2 \mid \forall X \in C_n(\operatorname{proj}\Lambda) \text{ indecomposable}, X^1 = 0 \text{ or } X^n = 0\}.$$

**Proof.** We noted the following:

i. If s.gl.dim($\Lambda$) = $\infty$, the result is clear.

ii. If s.gl.dim($\Lambda$) = $\eta$, then there is $X$ indecomposable complex in $C_{\eta+1}(\text{proj}\Lambda)$ such that $\ell(X) = \eta$. The complex $X$ satisfies that $X^1 \neq 0$ and $X^{\eta+1} \neq 0$; therefore,

$$\eta + 1 < m_0.$$

And by Lemma 3.1, we have $m_0 \leq \eta + 2$; therefore, $\eta = m_0 - 2$. $\square$

**Remark 3.23.** In the above theorem, it is not necessary that $\Lambda$ has finite strong global dimension.

Following to D. Vossieck in [22], $D^b(\Lambda)$ is called **discrete** and $\Lambda$ **derived-discrete** if for every vector $\boldsymbol{n} = (n_i)_{i \in \mathbb{Z}}$ of natural numbers there are only finitely many isomorphism classes of indecomposable objects in $D^b(\Lambda)$ of homological dimension vector $\boldsymbol{n}$.

If $\Lambda$ is derived discrete by R. Bautista in [[3], Theorem 2.3 (a)], then $C_n(\text{proj}\Lambda)$ is of finite representation type for all $n \geq 1$. Furthermore, if $\Lambda$ is piecewise hereditary by definition, there exists a hereditary abelian category $\mathcal{H}$ such that $D^b(\Lambda) \cong D^b(\mathcal{H})$ as triangulated categories. So, from D. Happel and D. Zacharia in [19], and from D. Vossieck in [22], if $\Lambda$ is piecewise hereditary and derived discrete, then $\mathcal{H}$ is given by a path algebra of Dynkin type, and we can calculate the strong global dimension of $\Lambda$ explicited using the following algorithm:

**Algorithm 3.24.** Let $\Lambda$ be a piecewise hereditary and derived discrete. In [12], C. Chaio, I, Pratti, and M.J. Souto Salorio showed how to build the Auslander-Reiten quiver of the category $C_n(\text{proj}\Lambda)$ for all $n > 1$.

Step 1. Build the Auslander-Reiten quiver of $C_2(\text{proj}\Lambda)$.

a) If all indecomposable complexes satisfy that the first or the second cells are zero, then by Theorem 3.22, we have s.gl.dim($\Lambda$) = $2 - 2 = 0$.

b) Otherwise, go to step 2.

Step 2. Build the Auslander-Reiten quiver of $C_3(\text{proj}\Lambda)$ and carry out the a) or b) of Step 1.

Step 3. Continue until obtaining the case (a). The algorithm stops because $\Lambda$ has finite strong global dimension by Theorem 2.2.

## 4. Applications for the derived category and examples

The aim of this section is to study some relations between $C_{\eta+1}(\text{proj}\Lambda)$ and the bounded derived category $D^b(\Lambda)$; we noted that if s.gl.dim($\Lambda$) = $\eta < \infty$, then $\Lambda$ has finite global dimension and so we can use the well-known fact $K^b(\text{proj}\Lambda) \cong D^b(\Lambda)$ as triangulated categories.

Let's recall from P. Grivel in [[17], Theorem 6.5], that $D^b(\Lambda)$ is a triangulated category. A triangle $X \xrightarrow{u} Y \xrightarrow{v} Z \xrightarrow{w} X[1]$ in $D^b(\Lambda)$ is called **Auslander-Reiten triangle** if the following conditions are satisfied:

(AR1) $X, Z$ are indecomposable.

(AR2) $w \neq 0$.

(AR3) If $f\colon W\to Z$ is not a retraction, then there is $f'\colon W\to Y$ such that $vf'=f$.

See D. Happel in [18], for a more detailed account on the Auslander-Reiten triangles.

Our first result allowed obtaining irreducible morphism and Auslander-Reiten triangles in $D^b(\Lambda)$ through the irreducible morphims and almost split $\mathcal{E}_{\eta+1}$-sequences in $C_{\eta+1}(\operatorname{proj}\Lambda)$ respectively.

Item i. of Theorem 4.1 was proved by C. Chaio, A.G. Chaio, and I. Pratti in [[10], Theorem 4.5], we will proof of this fact using the functors defined in section 2.

**Theorem 4.1.** *Let $\Lambda$ be an algebra with finite strong global dimension where $\eta = $ s.gl.dim$(\Lambda)$.*

  i. *Let $X$ and $Y$ be indecomposable complexes in $C_{\eta+1}(\operatorname{proj}\Lambda)$ without $\mathcal{E}_{\eta+1}$-projective and $\mathcal{E}_{\eta+1}$-injective direct summands. Suppose that $X$ can not be extended to the left in $C_{\eta+2}(\operatorname{proj}\Lambda)$ and $Y$ can not be extended to the right in $C_{\eta+2}(\operatorname{proj}\Lambda)$. Up to translation, $f\colon X\to Y$ is an irreducible morphism in $C_{\eta+1}(\operatorname{proj}\Lambda)$ if and only if $f\colon X\to Y$ is an irreducible morphism in $D^b(\Lambda)$.*

  ii. *Let $X,Y$ and $Z$ be complexes in $C_{\eta+1}(\operatorname{proj}\Lambda)$ without $\mathcal{E}_{\eta+1}$-projective and $\mathcal{E}_{\eta+1}$-injective direct summands. Given*

  $$X\xrightarrow{u}Y\xrightarrow{v}Z\xrightarrow{w}X[1] \tag{4.1}$$

  *a triangle in $D^b(\Lambda)$. Up to translation $X\xrightarrow{u}Y\xrightarrow{v}Z\xrightarrow{w}X[1]$ is Auslander-Reiten triangle in $D^b(\Lambda)$ if and only if there exists an almost split $\mathcal{E}_{\eta+1}$-sequence in $C_{\eta+1}(\operatorname{proj}\Lambda)$*

  $$X\xrightarrow{\left[\begin{smallmatrix}u\\u'\end{smallmatrix}\right]}Y\oplus P\xrightarrow{[\,v\ \ v'\,]}Z, \tag{4.2}$$

  *such that $P$ is an indecomposable complex $\mathcal{E}_{\eta+1}$-projective-injective and (4.2) induces an Auslander-Reiten triangle isomorphic to (4.1).*

**Proof.**

  i. Suppose that $f\colon X\to Y$ is an irreducible morphism in $C_{\eta+1}(\operatorname{proj}\Lambda)$. Since $X$ can not be extended to the left in $C_{\eta+2}(\operatorname{proj}\Lambda)$, $i^{\eta+1}_{-1}(f)$ is an irreducible morphism in $C_{\eta+2}(\operatorname{proj}\Lambda)$ by Proposition 2.7 *i*. We noted that $i^{\eta+1}_{-1}(Y)$ can not be extended to the right in $C_{\eta+3}(\operatorname{proj}\Lambda)$, $i^{\eta+2}_{0}(i^{\eta+1}_{-1}(f))$ is an irreducible morphism in $C_{\eta+3}(\operatorname{proj}\Lambda)$ by Proposition 2.7 *ii*. By R. Bautista and M.J. Souto Salorio in [[5], Corollary 4.5], we have that $i^{\eta+2}_{0}(i^{\eta+1}_{-1}(f))$ is an irreducible morphism in $D^b(\Lambda)$ and by Remark 2.3, up to translation, $f$ is irreducible in $D^b(\Lambda)$.

  Conversely, if $f\colon X\to Y$ is an irreducible morphism in $D^b(\Lambda)$, since $K^b(\operatorname{proj}\Lambda)\cong D^b(\Lambda)$ by H. Giraldo and H. Merklen in [[16], Theorem 6], (or also see [[5], Proposition 3.3 (b)]) $f$ is an irreducible morphism in $C_{[m,n]}(\operatorname{proj}\Lambda)$ for some $m,n\in\mathbb{Z}$ with $m<n$, by Remark 2.3, up to translation, $f$ is an irreducible morphism in $C_{n'}(\operatorname{proj}\Lambda)$ for some $n'\in\mathbb{Z}$ and finally, the statement follows from the Corollary 3.4 and Theorem 3.21 i.

  ii. Let $X,Y$ and $Z$ be complexes in $C_{\eta+1}(\operatorname{proj}\Lambda)$. We noted that, up to translation, $X\to Y\to Z$ is an almost split $\mathcal{E}_{\eta+1}$-sequence in $C_{\eta+1}(\operatorname{proj}\Lambda)$ if and only if $X\to Y\to Z$ is an almost split $\mathcal{E}_{[m,n]}$-sequence in $C_{[m,n]}(\operatorname{proj}\Lambda)$ for some $m,n\in\mathbb{Z}$ with $m<n$.

If $X \xrightarrow{\begin{bmatrix} u \\ u' \end{bmatrix}} Y \oplus P \xrightarrow{[\, v \; v' \,]} Z$ is an almost split $\mathcal{E}_{[m,n]}$-sequence in $\mathrm{C}_{[m,n]}(\mathrm{proj}\Lambda)$, then this induce the Auslander-Reiten triangle $X \xrightarrow{u} Y \xrightarrow{v} Z \to X[1]$ in $\mathrm{D}^b(\Lambda)$ by R. Bautista, M.J. Souto Salorio, and R. Zuazua in [[6], Lemma 9.1].

Conversely, let $X \xrightarrow{u} Y \xrightarrow{v} Z \xrightarrow{w} X[1]$ be an Auslander-Reiten triangle in $\mathrm{D}^b(\Lambda)$, then using [[8], Theorem 3.1], (or also see [[9], Theorem 2.7]) we have that there exists $X \xrightarrow{\begin{bmatrix} u \\ u' \end{bmatrix}} Y \oplus P \xrightarrow{[\, v \; v' \,]} Z$ an almost split $\mathcal{E}_{[m,n]}$-sequence in $\mathrm{C}_{[m,n]}(\mathrm{proj}\Lambda)$ for some $m, n \in \mathbb{Z}$ with $m < n$, where $P$ is an indecomposable complex $\mathcal{E}_{[m,n]}$-projective-injective, by Remark 2.3, up to translation, $X \xrightarrow{\begin{bmatrix} u \\ u' \end{bmatrix}} Y \oplus P \xrightarrow{[\, v \; v' \,]} Z$ is an almost split $\mathcal{E}_{n'}$-sequence in $\mathrm{C}_{n'}(\mathrm{proj}\Lambda)$ for some $n' \in \mathbb{Z}$ and finally, the statement follows from the Corollary 3.4 and Theorem 3.21 ii. □

We will define an special subquiver of Auslander-Reiten quiver of $\mathrm{C}_{\eta+1}(\mathrm{proj}\Lambda)$.

**Definition 4.2.** *Let $\Lambda$ be an algebra with finite strong global dimension where $\eta = \mathrm{s.gl.dim}(\Lambda)$ and $\Gamma_{\mathrm{C}_{\eta+1}(\mathrm{proj}\Lambda)}$ the Auslander-Reiten quiver of $\mathrm{C}_{\eta+1}(\mathrm{proj}\Lambda)$. We will define $\bar{\mathrm{L}}_{\mathrm{C}_{\eta+1}(\mathrm{proj}\Lambda)}$ by the quiver obtained from $\Gamma_{\mathrm{C}_{\eta+1}(\mathrm{proj}\Lambda)}$ in the following way, in each connected component of $\Gamma_{\mathrm{C}_{\eta+1}(\mathrm{proj}\Lambda)}$ we will remove the indecomposable complexes such that both are $\mathcal{E}_{\eta+1}$-projective and $\mathcal{E}_{\eta+1}$-injective and afterwards we will consider the connected component given by the indecomposable complex $X$ such that $X$ can not be extended neither to the left in $\mathrm{C}_{\eta+2}(\mathrm{proj}\Lambda)$ nor to the right in $\mathrm{C}_{\eta+2}(\mathrm{proj}\Lambda)$.*

We have the following theorem using the above definition:

**Theorem 4.3.** *If $\Lambda$ has finite strong global dimension with $\eta = \mathrm{s.gl.dim}(\Lambda)$, then the Auslander-Reiten quiver of $\mathrm{D}^b(\Lambda)$ is isomorphic to $\mathbb{Z}$ copies of $\bar{\mathrm{L}}_{\mathrm{C}_{\eta+1}(\mathrm{proj}\Lambda)}$.*

**Proof.** Up to translation and connecting map (in the sense of [[18], I.5.5]), we noted that $\Gamma_{\mathrm{D}^b(\Lambda)} \cong \mathbb{Z}\bar{\mathrm{L}}_{\mathrm{C}_{\eta+1}(\mathrm{proj}\Lambda)}$ by Theorem 4.1 ii, and Corollary 3.4. □

**Theorem 4.4.** *Let $\Lambda$ be an algebra with finite strong global dimension where $\eta = \mathrm{s.gl.dim}(\Lambda)$. $\mathrm{C}_{\eta+1}(\mathrm{proj}\Lambda)$ is of tame representation type if and only if $\mathrm{C}_n(\mathrm{proj}\Lambda)$ is of tame representation type for all $n \geq 1$.*

**Proof.** Suppose that $\mathrm{C}_{\eta+1}(\mathrm{proj}\Lambda)$ is of tame representation type.

i. If $n \leq \eta + 1$, by R. Bautista in [[3], Theorem 1.1], $\mathrm{C}_n(\mathrm{proj}\Lambda)$ is of tame representation type or wild representation type. If $\mathrm{C}_n(\mathrm{proj}\Lambda)$ is of wild representation type, as $\mathrm{C}_n(\mathrm{proj}\Lambda) \subseteq \mathrm{C}_{\eta+1}(\mathrm{proj}\Lambda)$, then $\mathrm{C}_{\eta+1}(\mathrm{proj}\Lambda)$ is of wild representation type, but this contradicts our hypothesis.

ii. If $n > \eta + 1$, we have the following fact:

$$\mathrm{indC}_n(\mathrm{proj}\Lambda) \subseteq \mathrm{indC}_{\eta+1}(\mathrm{proj}\Lambda) \bigcup \mathrm{indC}_{\eta+1}(\mathrm{proj}\Lambda)[1] \bigcup \cdots \bigcup \mathrm{indC}_{\eta+1}(\mathrm{proj}\Lambda)[n-1].$$

where $\mathrm{indC}_n(\mathrm{proj}\Lambda)$ denote the isomorphisms classes of the indecomposable complexes in $\mathrm{C}_n(\mathrm{proj}\Lambda)$ for all $n \geq 1$.

So the result is followed of Lemma 3.1 and our hypothesis.

Conversely, it is evident.   □

We complete this section given two examples in which we show how to calculate the strong global dimension of an algebra with the conditions given in Algorithm 3.24.

**Example 4.5.** We will calculate the strong global dimensional of the following finite dimensional $\Bbbk$-algebra. We will consider $\Lambda$ the path algebra given by the following quiver

$$Q \;:\; 1 \xrightarrow{\alpha} 2 \xrightarrow{\beta} 3$$

with $\alpha\beta = 0$.

By I. Reiten in [[21], Theorem 4.4], and D. Happel and D. Zacharia in [[19], Proposition 3.3], we know that s.gl.dim($\Lambda$) = 2. Now, we apply Algorithm 3.24. First, we built the Auslander-Reiten quiver of $C_2(\text{proj}\Lambda)$:

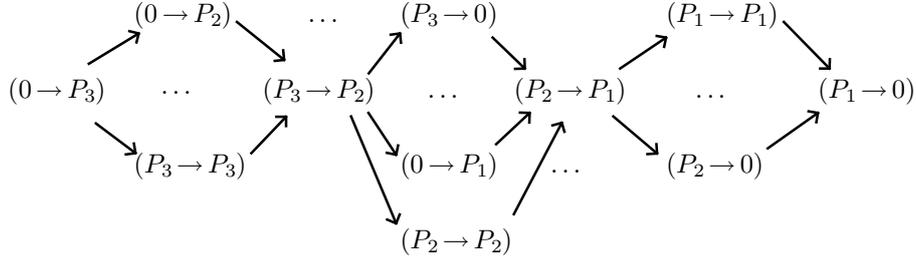

We noted that the complex $X := P_3 \longrightarrow P_2$ in $C_2(\text{proj}\Lambda)$ satisfies that the first and the last cells are non-zero. Therefore, we built the Auslander-Reiten quiver of $C_3(\text{proj}\Lambda)$:

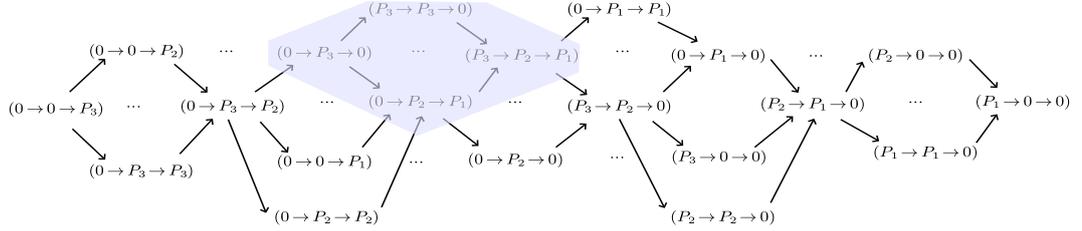

Now, the complex $X := P_3 \longrightarrow P_2 \longrightarrow P_1$ in $C_3(\text{proj}\Lambda)$ satisfies that the first and the last cells are non-zero, then we built the Auslander-Reiten quiver of $C_4(\text{proj}\Lambda)$

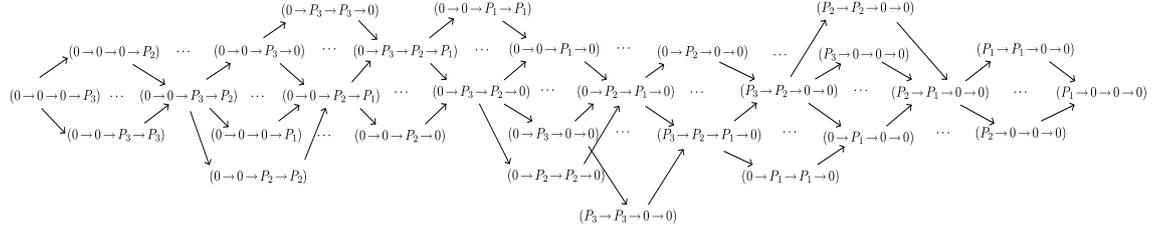

In this step, all complexes in $C_4(\text{proj}\Lambda)$ satisfy that the first or the last cells are zero, then by Theorem 3.22 we have s.gl.dim($\Lambda$) = 2.

**Remark 4.6.** In the example above

  i. Every irreducible morphism and every almost split $\mathcal{E}_4$-sequence of $C_4(\text{proj}\Lambda)$ are an irreducible morphism and an almost split $\mathcal{E}_3$-sequence of $C_3(\text{proj}\Lambda)$ and vice versa, up to translation. This is the result of Proposition 3.19.

ii. We noted that the shaded almost split $\mathcal{E}_3$-sequence in the above Auslander-Reiten quiver for $C_3(\text{proj}\Lambda)$ is not an almost split $\mathcal{E}_2$-sequence in the Auslander-Reiten quiver in $C_2(\text{proj}\Lambda)$, then we can say that there is one almost split $\mathcal{E}_3$-sequence in $C_3(\text{proj}\Lambda)$ such that this is not an almost split $\mathcal{E}_2$-sequence in $C_2(\text{proj}\Lambda)$. This fact can be generalized. There is one almost split $\mathcal{E}_{\eta+1}$-sequence in $C_{\eta+1}(\text{proj}\Lambda)$ such that it is not an almost split $\mathcal{E}_\eta$-sequence in $C_\eta(\text{proj}\Lambda)$. Since $\eta = \text{s.gl.dim}(\Lambda)$ then there exists $Z \in K^b(\text{proj}\Lambda)$ an indecomposable complex non-zero such that $\ell(Z) = \eta$. By the Remark 2.10 iii., $Z \in C_{\eta+1}(\text{proj}\Lambda)$ and $Z \notin C_\eta(\text{proj}\Lambda)$, again, by the Remark 2.10 i., $Z$ is not $\mathcal{E}_{\eta+1}$-indecomposable projective, then by R. Bautista, M.J. Souto Salorio, and R. Zuazua in [[6], Theorem 8.2] there is only an almost split $\mathcal{E}_{\eta+1}$-sequence in $C_{\eta+1}(\text{proj}\Lambda)$ that end in $Z$

$$X \to Y \to Z$$

Since $Z \notin C_\eta(\text{proj}\Lambda)$, then the above almost split $\mathcal{E}_{\eta+1}$-sequence is not in $C_\eta(\text{proj}\Lambda)$.

**Example 4.7.** In this example, we will calculate the strong global dimension of a finite dimensional $\Bbbk$-algebra using two tecniques. This example also shows that the strong global dimension is not the same as the global dimension. In [[7], Example 4.2.1], Y. Calderón-Henao calculated the s.gl.dim$(\Lambda)$ and gl.dim$(\Lambda)$ for $\Lambda$ the path algebra given by the following quiver

$$Q : 1 \xrightarrow{\alpha} 2 \xrightarrow{\beta} 3 \xrightarrow{\gamma} 4 \xrightarrow{\delta} 5 \xrightarrow{\eta} 6$$

with $\alpha\beta = \beta\gamma = \delta\eta = 0$.

That is,

$$\text{gl.dim}(\Lambda) = 3 \quad \text{and} \quad \text{s.gl.dim}(\Lambda) = 4.$$

Actually, Y. Calderón-Henao found the indecomposable complex with the longest length, which is

$$Z \colon P_6 \longrightarrow P_5 \longrightarrow P_3 \longrightarrow P_2 \longrightarrow P_1.$$

We observed that this tecnique is also used in [[11], Section 4] by C. Chaio, A.G. Chaio, and I. Pratti.

Now, we can calculate the strong global dimension using Algorithm 3.24. The Auslander-Reiten quivers of $C_n(\text{proj}\Lambda)$ for this algebra are extremely big, for this reason we do not draw them here. Nevertheless, we observed that the complex $Z$ is an indecomposable complex in $C_5(\text{proj}\Lambda)$ where the first and the last cells are non-zero. So we built the Auslander-Reiten quiver of $C_6(\text{proj}\Lambda)$ and we noted that in every indecomposable complex in $C_6(\text{proj}\Lambda)$ the first or the last cell is zero, so by Theorem 3.22 we also have that s.gl.dim$(\Lambda) = 4$.

## Acknowledgments

We thank CODI (Universidad de Antioquia, UdeA) and Fundación para la Promoción de la Investigación y la Técnología, entities which supported this research. We also thank Colciencias (Convocatoria 727 de 2015), which supported one of the authors of this research.

Instituto de Matemáticas, Universidad de Antioquia, Calle 67 No. 53-108, Medellín, Colombia. Instituto Tecnológico metropolitano-ITM
*Email address:* `yohny.calderon@udea.edu.co`
*Email address:* yohnycalderon8670@correo.itm.edu.co

Instituto de Matemáticas, Universidad de Antioquia, Calle 67 No. 53-108, Medellín, Colombia.
*Email address:* `cristianf.gallego@udea.edu.co`

Instituto de Matemáticas, Universidad de Antioquia, Calle 67 No. 53-108, Medellín, Colombia.
*Email address:* `hernan.giraldo@udea.edu.co`